\documentclass{amsart}

\usepackage{amsfonts}
\usepackage{amsmath}
\usepackage{amsthm}
\usepackage{amssymb}
\usepackage[english]{babel}
\usepackage{graphics}
\usepackage{latexsym}
\usepackage{longtable} 
\usepackage{mathrsfs} 
\usepackage{graphicx}
\usepackage{wrapfig} 
\usepackage[pdftex,colorlinks=true,linkcolor=blue,citecolor=cyan]{hyperref}

\newtheorem{theorem}{Theorem}
\newtheorem{corollary}[theorem]{Corollary}
\newtheorem{lemma}[theorem]{Lemma}

\newtheorem{claim}[theorem]{Claim}
\newtheorem{example}[theorem]{Example}

\theoremstyle{definition}
\newtheorem{definition}[theorem]{Definition}
\newtheorem{remark}[theorem]{Remark}

\newcommand{\mA}{\mathcal{A}}

\newcommand{\mF}{\mathcal{F}}

\newcommand{\mP}{\mathscr{P}}
\newcommand{\mM}{\mathcal{M}}
\newcommand{\mO}{\mathcal{O}}
\newcommand{\mY}{\mathscr{Y}}
\newcommand{\mD}{\mathcal{D}}

\newcommand{\E}{\mathrm{E}}
\newcommand{\D}{\mathrm{D}}
\renewcommand{\H}{\mathrm{H}}

\newcommand{\R}{\mathbb{R}}

\newcommand{\N}{\mathbb{N}}
\newcommand{\mS}{\mathbb{S}}
\newcommand{\mB}{\mathbb{B}}
\newcommand{\X}{{\bf X}}

\renewcommand{\P}{\mathrm{P}}

\newcommand{\noi}{\noindent}
\newcommand{\ms}{\medskip}

\newcommand{\al}{\alpha}
\newcommand{\be}{\beta}
\newcommand{\ga}{\gamma}

\newcommand{\de}{\delta}
\newcommand{\De}{\Delta}
\newcommand{\e}{\varepsilon}

\newcommand{\la}{\lambda}

\newcommand{\Om}{\Omega}

\newcommand{\weak }{\, -\!\!\!\!-\!\!\!\rightharpoonup}
\newcommand{\weakstar }{ \overset{\, *_{\phantom{|}}}{{\smash{\weak }}\, } }

\newcommand{\larrow}{\longrightarrow}
\newcommand{\ot}{\otimes}

\newcommand{\ri}{\rightarrow}

\newcommand{\p}{\partial}
\newcommand{\sub}{\subseteq}
\newcommand{\set}{\setminus}
\newcommand{\by}{\times}

\newcommand{\ess}{\mathrm{ess}}

\newcommand{\dist}{\mathrm{dist}}

\newcommand{\Div}{\mathrm{Div}}

\newcommand{\bt}{\begin{theorem}}\newcommand{\et}{\end{theorem}}
\newcommand{\bd}{\begin{definition}}\newcommand{\ed}{\end{definition}}
\newcommand{\bl}{\begin{lemma}}\newcommand{\el}{\end{lemma}}
\newcommand{\beq}{\begin{equation}}\newcommand{\eeq}{\end{equation}}
\newcommand{\bc}{\begin{claim}}\newcommand{\ec}{\end{claim}}
\newcommand{\bex}{\begin{example}}\newcommand{\eex}{\end{example}}
\newcommand{\bcor}{\begin{corollary}}\newcommand{\ecor}{\end{corollary}}
\newcommand{\bp}{\begin{proof}}\newcommand{\ep}{\end{proof}}

\newcommand{\BPL}{\medskip \noindent \textbf{Proof of Lemma} }

\newcommand{\BPCOR}{\medskip \noindent \textbf{Proof of Corollary} }

\newcommand{\BPT}{\medskip \noindent \textbf{Proof of Theorem} }

\numberwithin{equation}{section}

\begin{document}

\title[Characterisation of the PDE system of calculus of variation in $L^\infty$]{A Pointwise Characterisation of the PDE System of Vectorial Calculus of Variations in $L^\infty$}

\author{Birzhan Ayanbayev and Nikos Katzourakis}
\address{Department of Mathematics and Statistics, University of Reading, Whiteknights, PO Box 220, Reading RG6 6AX, Berkshire, UK}

  \thanks{\!\!\!\!\!\!\texttt{N.K. has been partially financially supported through the EPSRC grant EP/N017412/1}}

\email{n.katzourakis@reading.ac.uk \ (corresponding author)}

\email{b.ayanbayev@pgr.reading.ac.uk}


\subjclass[2010]{35D99, 35D40, 35J47, 35J47, 35J92, 35J70, 35J99}

\date{}


\keywords{$\infty$-Laplacian; Generalised solutions; Calculus of Variations in $L^\infty$; Young measures; fully nonlinear systems.}

\begin{abstract} Let $n,N\in \mathbb{N}$ with $\Om \subseteq \R^n$ open. Given $\H \in C^2(\Omega\times \mathbb{R}^N\times\R^{Nn}),$ we consider the functional 
\[ \tag{1} \label{1}
  \E_\infty (u,\mathcal{O})\, :=\, \underset{\mathcal{O}}{\mathrm{ess}\,\sup}\, \H (\cdot,u,\mathrm{D} u)  ,\ \ \ u\in W^{1,\infty}_\text{loc}(\Omega,\mathbb{R}^N),\ \ \ \mathcal{O} \Subset \Omega.
\]
The associated PDE system which plays the role of Euler-Lagrange equations in $L^\infty$ is
\[ \tag{2} \label{2}
\left\{
\begin{split}
\H_{P}(\cdot, u, \mathrm{D}u)\, \D\big(\H(\cdot, u, \mathrm{D} u)\big) \, = \, 0,\\ 
\H(\cdot, u, \mathrm{D} u) \, [\![\H_{P}(\cdot, u, \mathrm{D} u)]\!]^\bot \Big(\mathrm{Div}\big(\H_{P}(\cdot, u, \mathrm{D} u)\big)- \H_{\eta}(\cdot, u, \mathrm{D} u)\Big)\, =\, 0,
\end{split}
\right.
\]
where $[\![A]\!]^\bot := \mathrm{Proj}_{R(A)^\bot}$. Herein we establish that generalised solutions to \eqref{2} can be characterised as local minimisers of \eqref{1} for appropriate classes of affine variations of the energy. Generalised solutions to \eqref{2} are understood as $\mathcal{D}$-solutions, a general framework recently introduced by one of the authors.
\end{abstract}

\maketitle


\section{Introduction} \label{section1}

{\it Calculus of Variations} is the branch of Analysis which deals with the problem of finding and studying extrema of nonlinear functionals defined on certain infinite-dimensional topological vector spaces, as well as with describing these extrema through appropriate necessary and sufficient conditions. Such problems are called {\it variational} and are ubiquitous in nature, being also of paramount importance for other sciences. In most applications, the functional one wishes to study models some kind of ``energy" or ``action".

Let $\H \in C^2(\Om\times\R^N \times \R^{Nn})$ be a given function, where $\Om\sub \R^n$ is an open set and $n,N\in \N$. One of the most standard particular class of functionals of interest in Calculus of Variations has the form of 
\[
\ \ \E (u,\Om)\, :=\,  \int_\Om \H \big(x,u(x),\mathrm{D} u(x)\big) \, \mathrm{d}x
\]
defined on differentiable maps (i.e.\ vectorial functions) $u:\R^n \supseteq \Om\larrow \R^N$. In the above, $\R^{Nn}$ denotes the space of $N\by n$ matrices wherein the gradient matrix
\[
\D u(x)\, =\, \big( \D_i u_\al(x)\big)_{i=1,...,n}^{\al=1,...,N}\ \in \ \R^{Nn}
\]
of such maps is valued. We have also used the symbolisations $x=(x_1,..., x_n)^\top$, $u=(u_1,...,u_N)^\top$ and $\D_i\equiv \p/\p x_i$. Latin indices $i, j, k,...$ will run in $\{1,...,n\}$ and Greek indices $\al,\be,\ga,...$ will run in $\{1,...,N\}$, even if the range of summation is not explicitly mentioned. The simplest variational problem is to search for {\it minimisers} $u$ of $\E$, sought in a class $\mathscr{C}$ of differentiable maps $u$, subject to some kind of prescribed boundary condition on $\p\Om$ to avoid trivial minimisers. This means that any putative minimiser $u \in \mathscr{C}$, if it exists, should satisfy
\[
\E(u,\Om) \, \leq\, \E(v,\Om), \ \ \text{ for all $v \in \mathscr{C}$ with $u=v$ on }\p\Om.
\]
If such a minimiser exists, then the real function $t \mapsto \E(tv+(1-t)u)$ has a minimum at $t=0$ and should satisfy 
\[
\frac{\mathrm{d}}{\mathrm{d}t}\Big|_{t=0}\E\big(u+t(v-u)\big)\, =\, 0 .
\]
By the chain rule, this leads, at least {\it formally}, to the next necessary conditions, known as the {\it Euler-Lagrange system of Partial Differential Equations} (PDE):
\[
\sum_i \D_i \Big( \H_{P_{\al i}}(\cdot, u,\D u)\Big) \, =\, \H_{\eta_\al}(\cdot, u,\D u), \ \ \ \al=1,\ldots, N.
\]
In the above, the subscripts $\H_{P_{\al i}}, \H_{\eta_\al}$ denote the partial derivatives of $\H$ with respect to the respective variables $P_{\al i}$ and $\eta_\al$. Further, since the integral is additive with respect to the domain on which we integrate, it can be easily seen that if $u$ is a minimiser, then
\[
\E(u,\mO)\, \leq\, \E(v,\mO) , \ \ \text{ for all $v \in \mathscr{C}$ with $u=v$ on }\p\mO,
\]
where $\mO \Subset \Om$, namely $\overline \mO$ is a compact subset of $\Om$. The above weaker condition still suffices to derive the Euler-Lagrange system and any putative $u$ satisfying it is called an {\it absolute (or local) minimiser}.

The above discussion, although completely {\it formal}, nonetheless captures the quintessence of Calculus of Variations. However, one needs to use hardcore analytic tools to make rigorous the above formal reasoning. In particular, a central problem is that the minimisers are sought in a class of at most once differentiable maps, which the PDE is of second order and one has to devise a way to make sense of the PDE {\it weakly}, since second derivatives of $u$ may not exist! Such objects are called {\it generalised solutions}. Finding a efficient concept of generalised solution which allows one to prove that such a generalised object in fact exists and study its properties is a highly nontrivial part of the problem. A particular relevant question of great interest is to identify conditions on $\H$ allowing to {\it characterise variationally} the PDE system in terms of the functional, namely to provide sufficient as well as necessary conditions.

In this paper we are interested in the variational characterisation of the PDE system arising as the analogue of the Euler-Lagrange equations when one considers vectorial minimisation problems for {\it supremal functionals} of the form
\beq   \label{1.1}
\ \ \E_\infty (u,\mO)\, :=\,  \underset{x\in \mathcal{O}}{\mathrm{ess}\,\sup}\, \H \big(x,u(x),\mathrm{D} u(x)\big) ,\ \ \ \mO \Subset \Omega,
\eeq
defined on maps $u:\R^n \supseteq \Om\larrow \R^N$. This is in the spirit of the above discussion, but for the modern class of functionals as in \eqref{1.1}. The scalar case $N=1$ first arose in the work of G.\ Aronsson in the 1960s \cite{A1,A2} who initiated the area of {\it Calculus of Variations in the space $L^\infty$}. The field is fairly well-developed today and the relevant bibliography is vast. For a pedagogical introduction to the topic accessible to non-experts, we refer to \cite{K8}. 

The study of the vectorial case $N\geq 2$ started much more recently and the full system \eqref{1.2}-\eqref{1.4} first appeared in the paper \cite{K1} in the early 2010s and it is being studied quite systematically ever since (see \cite{K2}-\cite{K7}, \cite{K9}-\cite{K11}, as well as the joint works of the second author with Abugirda, Pryer, Croce and Pisante \cite{AK, CKP, KP, KP2}). The appropriate class of maps to place and study the functional is the Sobolev space $W^{1,\infty}(\Om,\R^N)$ of $L^\infty$ maps with $L^\infty$ derivative defined a.e.\ on $\Om$ (see e.g.\ \cite{E}). The direct extension of the concept of absolute minimisers for \eqref{1.1} reads
\[
\ \ \ \ \E_\infty(u,\mO)\,\leq\, \E_{\infty}(u+\phi,\mO), \ \ \ \ \mO \Subset \Om,\ \phi \in W^{1,\infty}_0(\mO,\R^N) 
\]
and was introduced and studied by Aronsson in the context of the scalar case. The subscript nought means that $\phi=0$ on $\p\mO$. The associated PDE system arising from  \eqref{1.1} as a necessary condition is
\beq \label{1.2}
\ \ \ \mF_{\infty}(\cdot, u,\D u,\D^2 u) \,=\, 0 \ \ \text{ in} \  \Om,
\eeq
where 
\[
\mF_\infty \ :\ \ \Om \by \R^N \by \R^{Nn} \by  \R^{Nn^2}_s \larrow \R^N
\]
is the Borel measurable map given by
\begin{align}
\label{1.3}
\mF_{\infty}(x,\eta,P,\X) \,&:=\, \H_{P}(x,\eta,P)\Big(\H_{P}(x,\eta,P):\X \, +
\,\H_{\eta}(x,\eta,P)^\top \!P\,+\,H_{x}(x,\eta,P) 
\Big) \nonumber\\
& \ +\,\H(x,\eta,P) \, [\![ \H_{P}(x,\eta,P) ]\!]^{\bot}\Big(\H_{PP}(x,\eta,P):\X\,+\,\H_{P\eta}(x,\eta,P) :P\\
& \ +\,\H_{Px}(x,\eta,P): \mathrm{I}\,  -\, \H_{\eta}(x,\eta,P) \Big).  \nonumber
\end{align}
In the above, $\R^{Nn^2}_s$ symbolises the space of symmetric tensors wherein the hessian of $u$ is valued:
\[
\D^2 u(x)\, =\, \left( \D^2_{ij}u_\al(x)\right)_{i,j =1,...,n}^{\al=1,...,N}\ \in \ \R^{Nn^2}_s.
\]
Further, $[\![A]\!]^\bot$ denotes the orthogonal projection onto the orthogonal complement of the range $R(A) \sub \R^N$ of a linear map $A :\R^n \larrow \R^N$:
\beq \label{1.4}
[\![A]\!]^\bot\, :=\, \mathrm{Proj}_{R(A)^\bot}.
\eeq
In index form, $\mF_\infty$ reads
\[
\begin{split}
\smash{\mF_{\infty}(x,\eta,\P,\X)_{\al}} &:= \sum_{i}\, \H_{P_{\al i}}(x,\eta,P) \Bigg(
\! \sum_{\be,j} \H_{P_{\be j}}(x,\eta,P)\X_{\be i j} + \sum_\be \H_{\eta_{\be}}(x,\eta,P)\, P_{\be i}
\\
&\ \ \ + \,  \H_{x_i}(x,\eta,P)\Bigg)\, + \ \H(x,\eta,P) \sum_{\be}\,  [\![\H_{P}(x,\eta,P) ]\!]^{\bot}_{\al\be} \, \cdot  
\\
& \ \ \ \cdot  \Bigg( \! \sum_{i,j} \H_{P_{\al i}P_{\be j}}(x,\eta,P) \, \X_{\be i j}
 \, + \, \sum_{i} \H_{P_{\al i}\eta_{\be}}(x,\eta,P)\, P_{\be i}
 \\ 
& \ \ \  + \, \sum_{i}  {\H_{P_{\al i}x_{i}}(x,\eta,P) \,-\, \H_{\eta_{\be}}(x,\eta,P) \Bigg)} ,
\end{split}
\]
where $\al =1,\ldots,N$. Note that, although $\H$ is $C^2$, the coefficient $[\![ \H_P(\cdot, u,\D u) ]\!]^\bot$ is \emph{discontinuous} at points where the rank of $\H_P(\cdot, u,\D u)$ changes. Further, because of the perpendicularity of $\H_P$ and $[\![ \H_P]\!]^\bot$ (that is $[\![ \H_P]\!]^\bot \H_P=0$), the system can be decoupled into the two independent systems
\[
\left\{\ \ 
\begin{split}
\H_{P}(\cdot, u, \mathrm{D}u) \, \D\big(\H(\cdot, u, \mathrm{D} u)\big) \, &= \, 0,
\\ 
\H(\cdot, u, \mathrm{D} u)\, [\![\H_{P}(\cdot, u, \mathrm{D} u)]\!]^\bot \Big(\mathrm{Div}\big(\H_{P}(\cdot, u, \mathrm{D} u)\big)- \H_{\eta}(\cdot, u, \mathrm{D} u)\Big)\, &=\, 0.
\end{split}
\right. \ \ \
\]
When $\H(x,\eta,P)=|P|^2$ (the Euclidean norm on $\R^{Nn}$ squared), the system \eqref{1.2}-\eqref{1.4} simplifies to the so-called {\it $\infty$-Laplacian}:
\beq \label{1.5}
\De_\infty u \, :=\, \Big(\D u \ot \D u + |\D u|^2[\![\D u]\!]^\bot \! \ot \mathrm{I} \Big):\D^2u\, =\, 0.
\eeq
In this paper we are interested in the characterisation of appropriately defined generalised vectorial solutions $u :\R^n \supseteq \Om \larrow \R^N$ to \eqref{1.2}-\eqref{1.4} in terms of the functional \eqref{1.1}. It is well known even from classical scalar considerations for $N=1$ that the solutions to \eqref{1.2}-\eqref{1.4} in general cannot be expected to be smooth. In the scalar case, generalised solutions are understood in the viscosity sense (see \cite{C,CIL,K8}). Since the viscosity theory does not work for \eqref{1.2}-\eqref{1.4} when $N\geq 2$, we will interpret solutions in the so-called {\it $\mD$-sense}. This is a new concept of generalised solutions for fully nonlinear systems of very general applicability recently introduced in \cite{K10}-\cite{K9}.

Deferring temporarily the details of this new theory of {\it $\mD$-solutions}, we stress the next purely vectorial peculiar occurrence: it is not yet known \emph{whether Aronsson's variational notion is appropriate when $\min\{n,N\}\geq 2$}. In the model case of \eqref{1.5} and for $C^2$ solutions, the relevant notion of so-called $\infty$-Minimal maps allowing to \emph{characterise variationally} solutions to \eqref{1.5} in term of $ u\mapsto \|\D u\|_{L^\infty(\cdot)}$ was introduced in \cite{K4}. These findings are compatible with the early vectorial observations made in \cite{BJW1, BJW2}, wherein the appropriate $L^\infty$ quasi-convexity notion in the vectorial case is essentially different from its scalar counterpart. In the recent paper \cite{K11} a new characterisation has been discovered that allows to connect $\mD$-solutions of \eqref{1.5} to local minimisers of $u\mapsto \|\D u\|_{L^\infty(\cdot)}$ in terms of certain classes of \emph{local affine variations}. This result offered new insights to the difficult problem of establishing connections of \eqref{1.1} to \eqref{1.2}-\eqref{1.4}.

In this paper we generalise the results of \cite{K11}, characterising general $\mD$-solutions to \eqref{1.2}-\eqref{1.4} in terms of local affine variations of \eqref{1.1}. Our main result is Theorem \ref{theorem11} that follows and asserts that $\mD$-solutions to \eqref{1.2}-\eqref{1.4} in $C^1(\Om,\R^N)$ can be characterised variationally in terms of \eqref{1.1}. The a priori $C^1$ regularity assumed for our putative solutions is slightly higher than the generic membership in the space $W^{1,\infty}(\Om,\R^N)$, but as a compensation we impose \emph{no convexity} of any kind for the hamiltonian $\H$ for the derivation of the system.

In special case of classical solutions, our result reduces to the following corollary which shows the geometric nature of our characterisation\footnote{We caution the reader that the statement of Corollary \ref{corollary1} sacrifices precision for the sake of clarity. The fully precise statement is that given in the main result, Theorem \ref{theorem11}.}:

\begin{corollary}[$C^2$ solutions of $\mF_{\infty}=0$]  \label{corollary1}  Let $\Om\sub \R^n$ be open, $u\in C^2(\Om,\R^N)$ and $\H \in C^2(\Om\times\R^n\times\R^{Nn})$. Then, 
\[
\ \ \ \mF_{\infty}(\cdot,u,\D u,\D^2u) \, = \, 0 \, \text{ in }\Om 
\ \ \Longleftrightarrow\ \ 
\left\{
\begin{array}{l}
\ \E_{\infty}(u,\,\mO)\, \leq\, \E_{\infty}(u+A,\,\mO)\, ,
\ms
\\
\ \forall\ \mO\Subset \Om, \ \forall\ A\in \big(\mA^{\parallel,\infty}_{\mO}\cup \mA^{\bot,\infty}_{\mO}\big)(u).
 \end{array}
 \right.
\]
Here $\mA^{\parallel,\infty}_{\mO}(u),\mA^{\bot,\infty}_{\mO}(u)$ are sets of affine maps given by
\[
\mA^{\parallel,\infty}_{\mO}(u)\, =\, \left\{A :\, \R^n \ri \R^N \, \left|
\begin{array}{l}
\D ^2A \equiv 0, \ A(x)=0 \text{ and exist } \xi\in\R^N \text{and}  \\ 
x\in \mO(u) \text{ s.t.\ the image of $A$ is parallel}\\
\text{to the tangent map of $\xi \H(\cdot,u,\D u)$ at }x 
 \end{array}
 \right. \!\!\!
 \right\}, \ \ \ 
\]
\[
\mA^{\bot,\infty}_{\mO}(u)\, =\, \left\{A :\, \R^n \ri \R^N \, \left|
\begin{array}{l}
\D^2A\, \equiv\, 0 \text{ and there exists $x\in \mO(u)$ s.t.\ the}   \\ 
\text{image of $A$ is normal to $\H_{P}(\cdot,u,\D u)$ at $x$ }\\
\text{and $A^\top \H_{P}(\cdot,u,\D u)$ is divergenceless at $x$ } 
 \end{array}
 \right. \!\!\!
 \right\}
\]
and 
\[
\mO(u)\,:=\, \mathrm{Argmax}\, \big\{\H(\cdot,u,\D u)\, :\, \overline{\mO} \big\}.
\]
\end{corollary}

This paper is organised as follows. In Section \ref{Section2} that follows we record all the basic facts needed regarding the concept of our $\mD$-solutions, namely our notion of generalised solution required to make rigorous sense of \eqref{1.2}-\eqref{1.4}. We also include a quick introduction to the analytic setup of so-called {\it Young measures}, on which $\mD$-solutions are based. We also give two simple auxiliary results which are utilised in the proof of our variational characterisation. Finally, in Section \ref{Section3} we state and prove our main result.

\section{Young measures, $\mD$-solutions and auxiliary results} \label{section2}
\label{Section2}

\noi \textbf{Young Measures.} Let $\Om\sub \R^n$ be open and $\mathbb{K}$ a compact subset of some Euclidean space $\R^d$. The set of {\it Young measures} $\mY\big(\Om, \mathbb{K} \big)$ forms a subset of the unit sphere of a certain $L^\infty$ space of measure-valued maps and this provides its useful properties, including sequential weak* compactness. More precisely, $\mY\big(\Om,  \mathbb{K} \big)$ is defined as
\[
\mY\big(\Om, \mathbb{K} \big) := \Big\{\nu :\, \Om \larrow \mathscr{P}( \mathbb{K} ) \, \Big| \ [\nu(\cdot)](\mathcal{U})\in L^\infty(\Om) \text{ for any open }\mathcal{U}\sub \mathbb{K} \Big\},
\]
where $\mathscr{P}( \mathbb{K} )$ is the set of probability measures on $\mathbb{K}$. To see how it arises, consider the separable  space $L^1 \smash{\big( \Om, C( \mathbb{K} ) \big)}$ of Bochner integrable maps. This space contains Carath\'eodory functions  $\Phi : \Om \by \mathbb{K} \larrow \R$ (namely functions for which $\Phi(\cdot,X)$ is measurable for all $X \in \mathbb{K}$ and $\Phi(x,\cdot)$ is continuous for a.e.\ $x\in\Om$) which satisfy 
\[
\| \Phi \|_{L^1 ( \Om, C( \mathbb{K} ))}:= \int_\Om  \big\| \Phi(x,\cdot)\big\|_{C^0(\mathbb{K} )} \mathrm{d}x \,<\, \infty. 
\]
We refer e.g.\ to \cite{FL, Ed, V} and to \cite{K9}-\cite{K11} for background material on these spaces. The dual space of this space is $L^\infty_{\mathrm{w}^*}\big( \Om, \mM( \mathbb{K} ) \big)$. This dual Banach space consists of Radon measure-valued maps $\smash{\Om \ni x \mapsto \nu(x)  \in \mM( \mathbb{K} )}$ which are weakly* measurable, in the sense that for any open set $\mathcal{U} \sub \mathbb{K}$, the function $x\mapsto [\nu(x)](\mathcal{U})$ is in $L^\infty(\Om)$. The norm of the space is given by 
\[
\| \nu \|_{ L^\infty_{\mathrm{w}^*} ( \Om, \mM ( \mathbb{K} ) ) }:=
\, \underset{x\in \Om}{\ess\,\sup} \left\|\nu(x) \right\|, 
\]
where ``$\|\cdot\|$" denotes the total variation. It thus follows that
\[
\mY\big(\Om, \mathbb{K} \big)\, =\, \Big\{ \nu\, \in \, L^\infty_{\mathrm{w}^*}\big( \Om,\mM( \mathbb{K} ) \big)\, : \, \nu (x) \in \mP ( \mathbb{K} ),\text{ for a.e. }x\in \Om\Big\}.
\]
\begin{remark}[Properties of Young Measures] \label{lemma0} We note the following facts about the set $\mY\big(\Om,\mathbb{K}\big)$ (proofs can be found e.g.\ in \cite{FG}): 

\smallskip
\noi i) It is convex and sequentially compact in the weak* topology induced from $L^\infty_{\mathrm{w}^*}$.  
\smallskip
 
\noi ii) The set of measurable maps $V : \R^n \supseteq \Om \larrow \mathbb{K}$ can be identified with a subset of it via the embedding $V \mapsto \de_V$, $\de_V(x):= \de_{V(x)}$. 
\smallskip

\noi iii)  Let $V^i,V^\infty : \R^n \supseteq \Om \larrow \mathbb{K}$ be measurable maps, $i\in \N$. Then, up the passage to subsequences, the following equivalence holds true as $i\ri\infty$: $V^{i} \larrow V^\infty$ a.e.\ on $\Om$ if and only if $\de_{V^i} \weakstar \de_{V^\infty}$ in $\mY\big(\Om, \mathbb{K}\big)$.
\end{remark}

\ms

\noi \textbf{$\mD$-solutions.} We now give some rudimentary facts about generalised solutions which are required for the main result in this paper. For simplicity we will restrict the discussion to $n=1$ for maps $u : \R \supseteq \Om \larrow \R^N$ with $\Om$ an interval. The notion of $\mD$-solutions is based on the probabilistic interpretation of limits of difference quotients by using \emph{Young} measures. Unlike standard PDE approaches which utilise Young measures valued in Euclidean spaces (see e.g.\ \cite{E, P, FL, CFV, FG, V, KR}), $\mD$-solutions are based on Young measures valued in the $1$-point compactification $\smash{\overline{\R}}^N := \R^N\cup \{\infty\}$ (which is isometric to the sphere $\mS^N$). The motivation of the notion in the case of $C^1$ solutions to $2$nd order fully nonlinear systems is the following: suppose temporarily $u\in C^2(\Om,\R^N)$ is a solution to 
\beq   \label{1.6}
\ \ \ \mF\big(x,u(x),u'(x),u''(x) \big)\,=\, 0, \quad \ x\in \Om,
\eeq
where $\mF : \Om \by \R^{N}\by \R^N\by \R^N\larrow\R^N$ is continuous. Let $\D^{1,h}$ be the usual difference quotient operator, i.e.\ $\smash{\D^{1,h} v(x) := \frac{1}{h}\big[ v(x+h) - v(x) \big]}$, $x\in\Om$, $h\neq0$. It follows that 
 \beq   \label{1.7}
\ \ \ \mF\Big(x,u(x),u'(x),\lim_{h\ri 0}\D^{1,h}u'(x) \Big)\,=\, 0, \quad \  x\in \Om.
\eeq
Since $\mF$ is continuous, \eqref{1.6} is equivalent to  
  \beq   \label{1.8}
\ \ \ \lim_{h\ri 0}\mF\Big(x,u(x),u'(x),\D^{1,h}u'(x) \Big)\,=\, 0, \quad \ x\in \Om.
\eeq
The crucial observation is that \emph{the limit in \eqref{1.8} may exist even if that of \eqref{1.7} does not}, whilst \eqref{1.8} makes sense for merely $C^1$ maps. In order to represent the limit in a convenient fashion, we need to view $u''$ and the difference quotients $\D^{1,h}u'$ as probability-valued maps  from $\Om$ to $\smash{\mathscr{P}\big( \smash{\overline{\R}}^N \big)}$, given by the respective Dirac masses $x \mapsto  \de_{\D^2u(x)}$ and $x \mapsto  \de_{\D^{1,h}u'(x)}$. The exact definition is as follows:

\begin{definition}[Diffuse Hessians]  \label{Diffuse hessians} Let $u :\R^n \supseteq \Om  \larrow\R^N$ be in $W^{1,\infty}(\Om,\R^N)$. Let also $D^{1,h}$ denote the difference quotient operator, i.e.\ $\smash{\D^{1,h} := \big(\D^{1,h}_1,...,\D^{1,h}_n \big)}$ and $\smash{\D^{1,h}_i v := \frac{1}{h}\big[ v(\cdot+he^i) - v \big]}$, $h\neq0$. The \textbf{diffuse hessians} $\mD^2 u$ of $u$ are the subsequential weak* limits of the difference quotients of the gradient in the set of sphere-valued Young measures along infinitesimal sequences $(h_{\nu})_{\nu=1}^\infty$:
\[
\ \ \ \ \ \de_{\D^{1,h_{\nu_k}} \D u} \weakstar \mD^2u \ \ \   \text{ in }\mY\big(\Om,  \smash{\overline{\R}}^{Nn^2}_s \big), \ \text{ as }k \ri \infty.
\]
\end{definition} 
Note that the set of Young measures is \emph{sequentially weakly*
compact} hence \emph{every map as above possesses diffuse $2$nd derivatives}.

\begin{definition}[$\mD$-solutions to $2$nd order systems] \label{definition13} Let $\Om\sub \R^n$ be an open set and $\mF :\Om \by \R^{N} \by \R^{Nn} \by \R^{Nn^2}_s  \larrow \R^N$ a Borel measurable map which is continuous with respect to the last argument. Consider the PDE system
\beq \label{2.11a}
\mF\big(\cdot,u,\D u,\D^2u\big)\, =\, 0\ \ \text{ on }\Om.
\eeq
We say that the locally Lipschitz continuous map $u : \R^n \supseteq \Om  \larrow \R^N$ is a \textbf{$\mD$-solution of \eqref{2.11a}} when for any diffuse hessian $\mD^2 u$ of $u$, we have
\beq \label{2.11aa}
\sup_{ \X_x  \in \, \mathrm{supp}_* (\mD^2 u(x)) }\,
\left| \mF\big(x,u(x),\D u(x), \X_x \big)\right| \, =\, 0, \quad \text{ a.e. }x\in\Om.
\eeq
Here ``$\mathrm{supp}_*$" symbolises the \emph{reduced support} of a probability measure excluding infinity, namely $\mathrm{supp}_* (\vartheta ) := \mathrm{supp}  ( \vartheta ) \set \{\infty\}$ when $\vartheta \in \mP\big( \smash{\overline{\R}}^{Nn^2}_s\big)$.
\end{definition}
We note that $\mD$-solutions are readily compatible with strong/classical solutions: indeed, by Remark \ref{lemma0}iii), if $u$ happens to be twice weakly differentiable then we have $\mD^2u(x)=\de_{\D^2u(x)}$ for a.e.\ $x\in\Om$ and the notion reduces to
\[
\ \ \ \sup_{\X_x\in \, \mathrm{supp}(\de_{\D^2 u(x)})} \big| \mF\big(x,u(x),\D u(x),\X_x\big)\big|\, =\, 0, \quad \text{ a.e. }x\in \Om,
\]
thus recovering strong/classical solutions because $\mathrm{supp}(\de_{\D^2 u(x)})=\{\D^2u(x)\}$. 

\ms

\noi \textbf{Two auxiliary lemmas.} We now identify two simple technical results which are needed for our main result.

\begin{lemma} \label{lemma1} Suppose $\Om\sub \R^n$ is open, $u\in C^1(\Om,\R^N)$  and $\H\in C^2(\R^n\times\R^N\times\R^{Nn})$. Fix $\mO\Subset \Om$ and an affine map $A :\R^n \larrow \R^N$. We set
\[
\mO(u)\, :=\, \Big\{ x\in \overline{\mO}\ : \ \H\big(x,u(x),\D u(x)\big) =\, \E_{\infty}(u,\,\mO)\Big\}.
\]
\noi a) If we have $\E_{\infty}(u,\,\mO) \leq \E_{\infty}(u+tA,\,\mO)$ for  all $t >0$, it follows that
\[
 \ \ \max_{z\in \overline{\mO}} \, \Big\{\H_{P}\big(z,u(z),\D u(z)\big) : \D A(z) \ +\  \H_{\eta}\big(z,u(z),\D u(z)\big) \cdot A(z) \Big\}\, \geq\, 0.
\]
In the above ``:" and ``$\cdot$" denote the inner products in $\R^{Nn}$ and $\R^N$ respectively.
\ms

\noi b) Let $x\in \mO$ and $0<\e <\dist(x,\p\mO)$. The set
\[
 \mO_\e (x)\, :=\, \Big\{ y\in \mO\, : \, \H(y,u(y),\D u(y))\leq \H\big(x,u(x),\D u(x)\big)\Big\}^{\!\circ} \bigcap \mB_\e(x)
\]
(where ``\,$(\cdot)^\circ$" denotes the interior) is open and compactly contained in $\mO$, whilst
\[
 \E_{\infty}\big(u,\, \mO_\e (x) \big)\, =\, \H\big(x,u(x),\D u(x)\big),
\]
whenever $\mO_\e (x) \neq \emptyset$.
\end{lemma}

\BPL \ref{lemma1}. a) Since $\E_{\infty}(u,\,\mO)\leq \E_{\infty}(u+tA,\,\mO)$, by Taylor-expanding $\H$, we have
\[
\begin{split}
0&\, \leq \,\max_{\overline{\mO}}\, \H\big(\cdot,u+tA, \D u \, + t\, \D A\big)  \, - \, \max_{\overline{\mO}}\, \H(\cdot,u,\D u) \\
&\, =\,\max_{\overline{\mO}}\,\Big\{\H(\cdot, u, \D u) \, +\, t\, \H_{\eta}(\cdot, u, \D u) \cdot A\, + \,t\, \H_{P}(\cdot, u, \D u) : \D A \\
&\, \ \ \ \ \, +\, O\big(t^2|A|^2  +\, t^2|\D A|^2\big)\Big\} \, -\,  \max_{\overline{\mO}}\, \H(\cdot,u,\D u)\\
&\, \leq \,  t\, \max_{\overline{\mO}} \,\Big\{ \H_{\eta}(\cdot, u, \D u)\cdot A \, +\, \smash{\H_{P}(\cdot, u, \D u) : \D A\Big\} \, +\, O(t^2) .}
\end{split}
\]
Consequently, by letting $t \ri 0$, we discover the desired inequality. Item b) is a direct consequence of the definitions.     \qed

\ms

Next, we have the following simple consequence of Danskin's theorem \cite{D}:

\begin{lemma} \label{lemma2} Given an open set $\Om\sub \R^n$, consider maps $u\in C^1(\Om,\R^N)$ and  $\H\in C^2(\R^n\times\R^N\times\R^{Nn})$, an affine map $A :\R^n \larrow \R^N$ and $\mO\Subset \Om$. We define
\[
\ \ \ r(\la)\, :=\, \E_{\infty}(u+\la A,\,\mO)\,-\,\E_{\infty}(u,\,\mO),\ \ \ \,\la\geq0.
\]
Let also $\mO(u)$ be as in Lemma \ref{lemma1}. Then, $r$ is convex, $r(0)=0$ and also it satisfies
\[
\ \ \ \underline{\D} r(0^+) \, \geq \ \max_{{\mO(u)}}\Big\{  \H_{P}(\cdot,u,\D u):\D A \ +\ \H_{\eta}(\cdot,u,\D u)\cdot A \Big\},
\]
where $\underline{\D} r(0^+):=\underset{\la\ri 0^+}{\lim\inf} \, \frac{r(\la)-r(0)}{\la} $ is the lower right Dini derivative of $r$ at zero.
\end{lemma}

\BPL \ref{lemma2}. The result is deducible from Danskin's theorem (see \cite{D}) but we prove it directly since the $1$-sided version above is not given explicitly in the paper. By setting
\[
R(\la,y)\, := \H\Big(y\, ,\, u(y)+\la A(y)\, , \, Du(y)+\,\la DA(y)\Big)
\]
we have $r(\la) = \max_{y\in \overline{\mO}}\, R(\la,y) - \max_{y\in \overline{\mO}} \,R(0,y)$, whilst for any $\la\geq0$ the maximum $\, \max_{y\in \overline{\mO}} \, R(\la,y)$ is realised at (at least one) point $y^{\la} \in \overline{\mO}$. Hence 
\[
\begin{split}
\frac{1}{\la}\big( r(\la)-r(0)\big)\, &=\, \frac{1}{\la}\Big[\max_{y\in \overline{\mO}}\, R(\la,y) \,-\, \max_{y\in \overline{\mO}} \,R(0,y)\Big] 
\\
 &=\, \frac{1}{\la}\Big[   R(\la,y^{\la}) \, -\, R(0,y^0) \Big] 
 \\
 & =\, \frac{1}{\la}\Big[ \big( R(\la,y^\la) -  R( \la,y^0) \big)\,+\,\big( R(\la,y^0)- R(0,y^0)\big)\Big]
\end{split} 
\]
and hence
\[
\begin{split}
\frac{1}{\la}\big( r(\la)-r(0)\big)\, &\geq \, \frac{1}{\la} \big( R(\la,y^0)\,-\, R(0,y^0) \big),
\end{split} 
\]
where $y^0\in \overline{\mO}$ is any point such that $R(0,y^0)\, =\, \max_{\overline{\mO}}R(0,\cdot)$. Hence, we have
\[
\begin{split}
\underline{\D} r(0^+)   &= \, \underset{\la\ri 0^+}{\lim\inf} \, \frac{1}{\la}\big( r(\la)-r(0)\big) \\
& \geq \, \max_{y^0 \in \overline{\mO}} \left\{  \underset{\la\ri 0^+}{\lim\inf} \, \frac{1}{\la} \Big( R(\la,y^0)- R(0,y^0) \Big) \right\} \\
& = \, \max_{y\in \mO(u)} \left\{ \underset{\la\ri 0^+}{\lim\inf} \, \frac{1}{\la} \Big( R(\la,y)- R(0,y) \Big) \right\}\\
& =\,  \max_{\mO(u)} \left\{ \underset{\la\ri 0^+}{\lim\inf} \, \frac{1}{\la} \bigg( \H\big(\, \cdot\, ,u + \la A,\D u +\la\D A\big) - \H\big(\cdot,u,\D u ) \bigg) \right\}\\
& = \, \max_{\mO(u)} \bigg\{ \underset{\la\ri 0^+}{\lim\inf} \, \frac{1}{\la} \bigg( \H (\cdot,u,\D u)\, +\, \la\, \H_{\eta}\big(\cdot,u,\D u\big)\cdot A \,  +\, \la \,\H_{P}( \cdot ,u,\D u) : \D A\\
&  \ \ \ \  + \, O\big(|\la \D A|^2+ |\la A|^2\big) -\, \H (\cdot,u,\D u) \bigg) \bigg\}
\end{split}
\]
and the desired inequality has been established. \qed

\smallskip

Let us  record the next simple inequality which follows from the definitions of lower right Dini derivative, in the case that $\H(x,\cdot,\cdot)$ is jointly convex for any $x\in\Om$. This is
\beq \label{2.1}
r(\la)\,-\,r(0)\,  \geq \, \underline{\D} r(0^+)\,\la, 
\eeq
for all $\la \geq 0$.

\section{The main result}
\label{Section3}

Now we proceed to the main theme of the paper, the variational characterisation of $\mD$-solutions to the PDE system \eqref{1.2} in terms of appropriate variations of the energy functional \eqref{1.1}. We recall that the Borel mapping $\mF_\infty : \Om \by \R^N \by \R^{Nn} \by \R^{Nn^2}_s \larrow \R^N$  is given by \eqref{1.3}-\eqref{1.4} and $\Om \sub \R^n$ is a fixed open set. 

\ms

\noi \textbf{Notational simplifications and perpendicularity considerations.} We begin by rewriting $\mF_{\infty}(\cdot ,u,\D u,\D^2 u)=0$ in a more malleable fashion. We define the maps
\begin{align}
 \mF_{\infty}^{\bot}(x,\eta,P,\X)\, &:= \,  \H_{PP}(x,\eta,P):\X\,+\,\H_{P\eta}(x,\eta,P): P \,+\,\H_{Px}(x,\eta,P):\mathrm{I} \phantom{\Big|} , 
    \label{3.1}
\\
\mF_{\infty}^{\parallel}(x,\eta,P,\X)\, &:=\,  \H_{P}(x,\eta,P):\X \,+\,\H_{\eta}(x,\eta,P)^\top P\,+\,H_{x}(x,\eta,P) 
    \label{3.2}
\end{align} 
and these are abbreviations of
\[
\begin{split}
\mF_{\infty}^{\bot}(x,\eta,P,\X)_\al\, &= \,  \sum_{\be, i, j}\H_{P_{\al i}P_{\be j}}(x,\eta,P)\, \X_{\be i j}\,+\,\sum_{\be , i} \H_{P_{\al i}\eta_\be}(x,\eta,P) \, P_{\be i} 
\\
&\ \ \ +\,\sum_i \H_{P_{\al i}x_i}(x,\eta,P) \phantom{\Big|} ,  
\\
\mF_{\infty}^{\parallel}(x,\eta,P,\X)_i\, &=\, \sum_{\be, j} \H_{P_{\be j}}(x,\eta,P)\X_{\be i j } \,+\,\sum_{\be }\H_{\eta_\be }(x,\eta,P)\, P_{\be i}\,+\,H_{x_i}(x,\eta,P). 
\end{split}
\]
Note that $\mF_{\infty}^{\bot}(x,\eta,P,\X) \in \R^N$, whilst $\mF_{\infty}^{\parallel}(x,\eta,P,\X) \in \R^n$. By utilising \eqref{3.1}-\eqref{3.2}, we can now express \eqref{1.3} as 
\[
\begin{split}
\mF_{\infty}(x,\eta,P,\X) \,&:= \, \H_{P}(x,\eta,P)\, \mF_{\infty}^{\parallel}(x,\eta,P,\X) + \, \H(x,\eta,P)\, \cdot
\\
& \ \ \ \ \cdot\,  [\![\H_{P}(x,\eta,P) ]\!]^{\bot}\Big(\mF_{\infty}^{\bot}(x,\eta,P,\X)\, -\, \H_{\eta}(x,\eta,P)\Big)   .
\end{split}
\]
Further, recall that in view of \eqref{1.4}, $[\![\H_{P}(x,\eta,P)]\!]^\bot$ is the projection on the orthogonal complement of $R(\H_{P}(x,\eta,P))$. Hence, by the orthogonality of $[\![\H_{P}(x,\eta,\!P)]\!]^{\bot}\cdot$    $\cdot\big(\mF_{\infty}^{\bot}(x,\eta,\!P,\X)-\H_{\eta}(x,\eta,\!P)\big)$ and $\H_{P}(x,\eta,P)\mF_{\infty}^{\parallel}(x,\eta,\!P,\X)$, we have 
\[
\ \ \ \ \ \ \ \ \ \mF_{\infty}(x,\eta,P,\X)= 0, \ \text{ for some }(x,\eta,P,\X) \in \,\Om \by \R^N \! \by \R^{Nn} \!\by \R^{Nn^2}_s ,
\]
if and only if
\[
\ \ \left\{\ \ 
\begin{split}
\H_{P}(x,\eta,P)\, \mF_{\infty}^{\parallel}(x,\eta,P,\X) \, &=\, 0,\\
\!\!\H(x,\eta,P) \, [\![\H_{P}(x,\eta,P) ]\!]^{\bot}\Big(\mF_{\infty}^{\bot}(x,\eta,P,\X)\!-\!\H_{\eta}(x,\eta,P)\Big)  \,&=\,0.
\end{split}
\right.
\]
Finally, for the sake of clarity we state and prove our characterisation below only in the case of $C^1$ solutions, but due to its pointwise nature, the result holds true for piecewise $C^1$ solutions with obvious adaptations which we refrain from providing. We will assume that the Hamiltonian $\H$ satisfies 
\beq \label{assH}
\ \ \ \ \ \ \smash{\big\{\H_P(x,\eta,\cdot) =0 \big\} \, \sub \, \big\{\H(x,\eta,\cdot) =0 \big\}, \ \ (x,\eta) \in \Om \by \R^N}.
\eeq
We will also suppose that the next set has vanishing measure
\beq \label{assm}
\Big|\Big\{x\in\Om\, :\, \mB_{r_x}(x)\bigcap \big\{h>h(x)\big\} \text{ is dense in }\mB_{r_x}(x)  \Big\}\Big|\, =\, 0,
\eeq
where $r_x \equiv {\dist(x,\p\Om)}$ and $h\equiv \H(\cdot,u,\D u)$. This assumption is natural, in the sense that it is satisfied by all know examples of explicit solutions. It is trivially satisfied if $h$ has no strict local minima in the domain.

\ms

Our main result is as follows:

\begin{theorem}[Variational characterisation of the PDE system arising in $L^\infty$] \label{theorem11} Let $\Om\sub \R^n$ be open, $u\in C^1(\Om,\R^N)$ and $\H \in C^2(\Om\times\R^n\times\R^{Nn})$ a function satisfying \eqref{assH} and suppose that \eqref{assm} holds. Then:

\noi $\mathrm{(A)}$ We have
\[
\phantom{\Big|} \mF_{\infty}(\cdot,u,\D u,\D^2u) =0 \ \text{ in }\Om,
\]
in the $\mD$-sense, if and only if 
\[
\ \ \ \ \phantom{\Big|}\E_{\infty}(u,\,\mO)\, \leq\, \E_{\infty}(u+A,\,\mO), \ \ \ \ \forall\  \mO\Subset \Om,\ \forall\ A\in \mA^{\parallel,\infty}_{\mO}(u) \bigcup \mA^{\bot,\infty}_{\mO}(u).
\]
For the sufficiency of the PDE for the variational problem we require that $\H(x,\cdot,\cdot)$ be convex. In the above, the sets $\mA^{\parallel,\infty}_{\mO}(u),\mA^{\bot,\infty}_{\mO}(u)$ consist, for any $\mO\Subset\Om$, by affine mappings as follows:
\[
\mA^{\parallel,\infty}_{\mO}(u) := \! \left\{ \! A :\, \R^n \ri \R^N  
\left|
\begin{array}{l}
\! \D ^2A \equiv0,\, A(x)=0\, \& \, \text{exist } \xi\in \R^N \!,\,  x\in \mO(u),\, \\ 
\! \mD^2u \in \mY\big(\Om,   \smash{\overline{\R}}^{Nn^2}_s \big)\ \&\ \X_x  \in \mathrm{supp}_*\big(\mD^2u(x)\big) \\
\text{s.t.  : }  \D A\, \equiv \,\xi \ot \mF_{\infty}^{\parallel}\big( x,u(x),\D u(x),\X_x \big)
 \end{array}
 \right. \!\!\!\!\!
 \right\} \! \bigcup  \R^N
\]
and
\[
\mA^{\bot,\infty}_{\mO}(u) := \!\left\{ \! A :\, \R^n \ri \R^N  
\left|  
\begin{array}{l}
\D^2A\, \equiv\, 0\ \& \text{ there exist } x\in \mO(u),\, \mD^2u\ \\ 
\in \mY\big(\Om,   \smash{\overline{\R}}^{Nn^2}_s \big)\ \&\, \X_x  \in \mathrm{supp}_*\big(\mD^2u(x)\big)
\\
\text{s.t.  : } A(x) \in R\Big(\H_{P}\big(x,u(x),\D u(x)\big)\Big)^\bot\ \\
\& \ \D A \in \mathscr{L} \big( x,  A(x) , \X_x \big)
 \end{array}
 \right. \!\!\!\!\!
 \right\}\!  \bigcup  \R^N
\]
where $\mathscr{L} \big( x,\eta, \X  \big)$ is an affine space of $N\by n$ matrices, defined as
\[
\mathscr{L}\big( x,\eta, \X \big) := 
\left\{
\begin{array}{l}
\!\!\Big\{ Q \in \R^{Nn}\ \Big| \ \, \H_{P}\big(x,u(x),\D u(x)\big)\!: Q   \ms\\
=\, -\eta \cdot \mF_{\infty}^{\bot}\big(x,u(x),\D u(x),\X\big)\Big\} \, ,
\ \,\ \  \,   \text{if}\ \H_{P}\big(x,u(x),\D u(x)\big) \neq \, 0, 
\ms
\\
\{0\}, \hspace{134pt} \, \text{ if}\ \H_{P}\big(x,u(x),\D u(x)\big) =\, 0,
\end{array}
\right.
\]
for any $(x,\eta,\X) \in  \Om \by \R^N \by \R^{Nn^2}_s$.

\smallskip

\noi $\mathrm{(B)}$ In view of the mutual perpendicularity of the two components of $\mF_\infty$ (see \eqref{3.1}-\eqref{3.2}), $\mathrm{(A)}$ is a consequence of the following particular results:
\[
\ \ \ \ \phantom{\Big|} \H_{P}(\cdot,u,\D u)\, \mF_{\infty}^{\parallel}(\cdot,u,\D u,\D^2u) \,=\,0  \ \text{ in }\Om, 
\]
in the $\mD$-sense, if and only if 
\[
\ \ \ \ \ \ \E_{\infty}(u,\,\mO)\, \leq\, \E_{\infty}(u+A,\,\mO), \ \ \ \forall\ \mO\Subset \Om,\  \forall\ A\in \mA^{\parallel,\infty}_{\mO}(u)
\]
and also
\[
\ \ \ \ \phantom{\Big|} \H(\cdot,u,\D u)\, [\![\H_{P}(\cdot,u,\D u)]\!]^{\bot}\Big(\mF_{\infty}^{\bot}(\cdot,u,\D u,\D^2u)-\H_{\eta}(\cdot,u,\D u)\Big)\,=\, 0 \ \text{ in }\Om,
\]
in the $\mD$-sense, if and only if 
\[
\ \ \ \ \ \ \E_{\infty}(u,\,\mO)\, \leq\, \E_{\infty}(u+A,\,\mO), \ \ \  \forall\ \mO\Subset \Om,\  \forall\ A\in \mA^{\bot,\infty}_{\mO}(u).
\]
\end{theorem}

We note that in the special case of $C^2$ solutions, Corollary \ref{corollary1} describes the way that classical solutions $u :\R^n \supseteq \Om \larrow \R^N$ to \eqref{1.2}-\eqref{1.4} are characterised.

\begin{remark}[About pointwise properties of $C^1$ $\mD$-solutions] \label{Remark} Let $u:\R^n \supseteq \Om \larrow \R^N$ be a $\mD$-solution to \eqref{1.2}-\eqref{1.4} in $C^1(\Om,\R^N)$. By Definition \ref{definition13}, this means that for any $\smash{ \mD^2u \in \mY\big(\Om,   \smash{\overline{\R}}^{Nn^2}_s \big)}$,
\[
\ \ \ \ \ \ \ \ \mF_\infty\big(x,u(x),\D u(x),\X_x\big)\, =\, 0, \ \ \text{ a.e. }x\in\Om, \ \forall\ \X_x \in \mathrm{supp}_*\big(\mD^2u(x)\big).
\]
By Definition \ref{Diffuse hessians}, every diffuse hessian of a putative solution is defined a.e.\ on $\Om$ as a weakly* measurable probability valued map $\R^n \supseteq \Om \larrow \mP \big( \smash{{\R}}^{Nn^2}_s\cup\{\infty\}\big)$. Let $\Om \ni x \mapsto  \textbf{O}_x  \in \smash{\R^{Nn^2}_s}$ be any selection of elements of the zero level sets
\[
\Big\{\X \in \smash{\R^{Nn^2}_s}\ : \ \mF_\infty\big(x,u(x),\D u(x), \X \big) =\, 0 \Big\}.
\]
By modifying each diffuse hessian on a Lebesgue nullset and choosing the representative which is redefined as $\mD^2u(x) = \de_{\textbf{O}_x}$ for a negligible set of $x$'s, we may assume that $\mD^2u(x)$ exists for all $x\in \Om$. Further, given that $\D u(x)$ exists for all $x\in \Om$, by perhaps a further re-definition on a Lebesgue nullset, it follows that $u$ is $\mD$-solution to \eqref{1.2}-\eqref{1.4} if and only if for (any such representative of) any diffuse hessian 
\[
\ \ \ \ \ \ \ \ \mF_\infty\big(x,u(x),\D u(x),\X_x\big)\, =\, 0, \ \ \ \ \forall\ x\in\Om,\  \forall\ \X_x \in \mathrm{supp}_*\big(\mD^2u(x)\big).
\] 
Note that at points $x\in \Om$ for which  $\mD^2u(x) =\de_{\{\infty\}}$ and hence $\mathrm{supp}_*\big(\mD^2u(x)\big)$ $= \emptyset$, the solution criterion is understood as being trivially satisfied. 
\end{remark}

\BPT \ref{theorem11}. It suffices to establish only (B), since (A) is a consequence of it. Suppose that for any $\mO\Subset \Om$ and any $A\in \smash{\mA^{\bot,\infty}_{\mO}(u)}$ we have $E_{\infty}(u,\mO) \leq \E_{\infty}(u+A,\mO)$. Fix a diffuse hessian $\mD^2u \in \mY\big(\Om, \smash{\overline{\R}}^{Nn^2}_s \big)$, a point $x\in \overline{\mO}$ such that $\mathrm{supp}_*\big(\mD^2u(x)\big) \neq \emptyset$ and an $\X_x \in \mathrm{supp}_*\big(\mD^2u(x)\big)$. In view of \eqref{3.1}, if $\H_{P}\big(x,u(x),\D u(x)\big)=0$, then, by our assumption on the level sets of $\H$, we have $\H\big(x,u(x),\D u(x)\big)=0$ as well and as a consequence we readily obtain
\beq \label{3.3}
\begin{split}
& \H\big(x,u(x),\D u(x)\big) \,[\![\H_{P}\big(x,u(x),\D u(x)\big)]\!]^\bot \cdot
\\
& \cdot 
  \Big(\mF^{\bot}_{\infty}\big( x,u(x),\D u(x),\X_x \big)-\H_{\eta}\big(x,u(x),\D u(x)\big)\Big)=0 
\end{split}
\eeq
is clearly satisfied at $x$. If $\H_{P}\big(x,u(x),\D u(x)\big)\neq 0$, then we select any direction normal to the range of $\H_{P}\big(x,u(x),\D u(x)\big) \in \R^{Nn}$, that is 
\[
n_x\, \in R\Big(\H_{P}\big(x,u(x),\D u(x)\big)\Big)^\bot \sub \,\R^N
\]
which means $n_x^\top \H_{P}\big(x,u(x),\D u(x)\big)=0$. Of course it may happen that the linear map $\H_{P}\big(x,u(x),\D u(x)\big) : \R^n \larrow \R^{Nn}$ is surjective and then only the trivial $n_x=0$ exists. In such an event, the equality \eqref{3.3} above is satisfied at $x$ because $[\![\H_{P}\big(x,u(x),\D u(x)\big)]\!]^\bot =0$. Hence, we may assume $n_x\neq 0$. Further, fix any matrix $N_x$ in the affine space $\mathscr{L}(x,n_x,\X_x ) \sub \R^{Nn}$. By the definition of $\mathscr{L}(x,n_x,\X_x)$, we have
\[
 \H_{P}\big(x,u(x),\D u(x)\big) : N_x\, =\, - n_x \cdot \mF_{\infty}^{\bot}\big( x,u(x),\D u(x),\X_x \big).
\]
Consider the affine map defined by
\[
\ \ \ A(z)\, :=\, n_x\, +\, N_x (z-x),\ \ \ z\in \R^n.
\]
We remark that $t A\in \mA^{\bot,\infty}_{\mO}(u)$ for any $t \in\R$. Indeed, this is a consequence of our choices and the next homogeneity property of the space $\mathscr{L}(x,\eta,\X)$: 
\[
\mathscr{L}(x, t \eta, \X) \, =\, t\, \mathscr{L}(x,\eta,\X), \ \ \ t \in \R. 
\] 
Let $\e>0$ be small, fix $x\in\Om$ and let us choose as $\mO$ the domain $\mO_\e(x)$ defined in Lemma \ref{lemma1}b). Our assumption \eqref{assm} implies that $\mO_\e(x)\neq \emptyset$ for a.e.\ $x\in\Om$. In view of the above considerations, we have
\[
\E_{\infty}\big(u,\mO_\e(x)\big) \, \leq \, \E_{\infty}\big(u+tA,\mO_\e(x)\big). 
\]
By applying Lemma \ref{lemma1}a), we have
\[
\begin{split}
0 \, &\leq\ \max_{z\in \overline{\mO_\e(x)}} \Big\{ \H_{P}\big( z,u(z),\D u(z) \big) : \D A(z) \, +\, \H_{\eta}\big( z,u(z),\D u(z) \big)\cdot A(z)\Big\}  
\\
& \!\!\!\!\! \overset{\e \ri 0}{-\!\!-\!\!\!\larrow} \ \H_{P}\big(x,u(x),\D u(x)\big) : N_x \, +\, \H_{\eta}\big(x,u(x),\D u(x)\big) \cdot n_x 
\\
& =\ -\, n_x \cdot \Big(\mF_{\infty}^{\bot}\big( x,u(x),\D u(x),\X_x \big)-\H_{\eta}\big(x,u(x),\D u(x)\big)\Big).
\end{split}
\] 
As a result, we have 
\[
n_x \cdot \Big(\mF_{\infty}^{\bot}\big( x,u(x),\D u(x),\X_x \big)-\H_{\eta}\big(x,u(x),\D u(x)\big)\Big)\leq 0
\]
for any direction $n_x \, \bot\, R\Big(\H_{P}\big(x,u(x),\D u(x)\big)\Big)$ and by the arbitrariness of $n_x$, we deduce that
\[
[\![\H_{P}\big(x,u(x),\D u(x)\big)]\!]^\bot\Big(\mF_{\infty}^{\bot}\big( x,u(x),\D u(x),\X_x \big) \,- \, \H_{\eta}\big(x,u(x),\D u(x)\big)\Big) \, = \, 0,
\]
for any $\smash{\mD^2 u \in \mY\big( \Om, \smash{\overline{\R}}^{Nn^2}_s \big)}$, $x\in \Om$ and $\X_x \in \mathrm{supp}_*\big(\mD^2 u(x)\big)$, as desired.

\ms

For the tangential component of the system we argue similarly.  Suppose that for any $\mO\Subset \Om$ and any $\smash{A\in \mA^{\parallel,\infty}_{\mO}(u)}$ we have $\E_{\infty}(u,\,\mO) \leq \E_{\infty}(u+A,\,\mO)$. Fix $x\in \overline{\mO}$, a diffuse hessian $\mD^2u \in \mY\big(\Om,   \smash{\overline{\R}}^{Nn^2}_s \big)$ such that $\mathrm{supp}_*\big(\mD^2u(x)\big)$ 
$\neq \emptyset$, a point $\X_x \in \mathrm{supp}_*\big(\mD^2u(x)\big)$ and $\xi \in \R^N$. Recalling \eqref{3.2}, we define the affine map
\[
\ \ \ A(z)\, :=\ \xi \ot \mF_{\infty}^{\parallel}\big(x,u(x),\D u(x),\X_x\big) \cdot (z-x),\ \ \ z\in \R^n. \\
\]
Fix $\e>0$ small, $x\in\Om$ and choose as $\mO$ the domain $\mO_\e(x)$ of Lemma \ref{lemma1}b). Then,  $t A \in \mA^{\parallel,\infty}_{\mO_\e(x)}(u)$ for any $t \in \R$. Consequently, in view our the above we have
\[
\E_{\infty}\big(u,\mO_\e(x)\big) \, \leq \, \E_{\infty}\big(u+tA,\mO_\e(x)\big) 
\]
and by applying Lemma \ref{lemma1}a), this yields
\[
\begin{split}
0 \, &\leq \ \max_{z\in \overline{\mO_\e(x)}} \Big\{\, \H_{P}\big( z,u(z),\D u(z) \big): \D A(z) \, + \, \H_\eta \big( z,u(z),\D u(z) \big) \cdot A(z) \Big\} \\
& \!\!\!\!\! \overset{\e \ri 0}{-\!\!-\!\!\!\larrow} \  \H_{P}\big(x,u(x),\D u(x)\big) : \Big(\xi \ot \mF_{\infty}^{\parallel}\big( x,u(x),\D u(x),\X_x \big)\Big).
\end{split}
\] 
Hence,
\[
 \xi \cdot \Big( \H_{P}\big(x,u(x),\D u(x)\big) \, \mF_{\infty}^{\parallel}\big( x,u(x),\D u(x),\X_x \big)\Big)\, \geq \, 0,
\] 
for any $\xi \in \R^N$. By the arbitrariness of $\xi$ we infer that 
\[
\H_{P}\big(x,u(x),\D u(x)\big) \, \mF_{\infty}^{\parallel}\big( x,u(x),\D u(x),\X_x \big)=0
\]
for any $\smash{\mD^2 u \in \mY\big(\Om,  \smash{\overline{\R}}^{Nn^2}_s \big)}$, $x\in \Om$ and $\X_x \in \mathrm{supp}_*\big(\mD^2 u(x)\big)$, as desired.

\smallskip

Conversely, let us fix $\mO\Subset \Om$, $x\in \mO(u)$, $\mD^2u \in \mY\big(\Om, \smash{\overline{\R}}^{Nn^2}_s \big)$, $\X_x \in \mathrm{supp}_*(\mD^2u(x))$ and $\xi \in \R^N$ corresponding to a map $A\in \smash{\mA^{\parallel,\infty}_{\mO}(u)}$. Let $r$ be the function of Lemma \ref{lemma2}. By applying Lemma \ref{lemma2} to the above setting, we have
\[
\begin{split}
\underline{\D}r(0^+)\, &\geq\, \max_{y\in {\mO(u)}} \, \Big\{ \,\H_{P}(y,u(y),\D u(y)) : \D A(y)  \ +\ \H_{\eta}(y,u(y),\D u(y))\cdot A(y)\Big\} \\
& \geq \, \,\H_{P}\big(x,u(x),\D u(x)\big) : \D A(x) \ +\  \H_{\eta}\big(x,u(x),\D u(x)\big)\cdot A(x)
\\
&=\ \H_{P}\big(x,u(x),\D u(x)\big):\Big(\xi \ot \mF_{\infty}^{\parallel}\big( x,u(x),\D u(x),\X_x \big)\Big)\\ 
&=\  \xi \cdot\Big( \H_{P}\big(x,u(x),\D u(x)\big) \mF_{\infty}^{\parallel}\big( x,u(x),\D u(x),\X_x \big)\Big)
\end{split}
\]
and hence $\underline{\D}r(0^+) \geq 0$ because $u$ is a $\mD$-solution. Due to the fact that $r(0)=0$ and $r$ is convex, by inequality \eqref{2.1} we have $r(t)\geq 0$ for all $t \geq 0$. Therefore,
\[
\ \ \ \E_{\infty}(u,\,\mO)\, \leq\, \E_{\infty}(u+A,\,\mO), \ \ \ \forall\ \mO\Subset\Om, \ \forall\ A\in \mA^{\parallel,\infty}_{\mO}(u).
\]
The case of $A\in \mA^{\bot,\infty}_{\mO}$ is completely analogous. Fix $\mD^2u \in \mY\big(\Om, \smash{\overline{\R}}^{Nn^2}_s \big)$, $\mO \Subset \Om$, $x\in \mO(u)$, $\X_x \in \mathrm{supp}_*(\mD^2u(x))$ and an $A$ with $A(x) \, \bot\, R\big(\H_{P}\big(x,u(x),\D u(x)\big)\big)$ and $\D A \in \mathscr{L}\big( x, A(x), \X_x\big)$. By applying Lemma \ref{lemma2} again, we have
\[
\begin{split}
\underline{\D}r(0^+)\, &\geq\ \max_{y\in {\mO(u)}} \, \Big\{\H_{P}(y,u(y),\D u(y)) : \D A(y)\ +\  \H_{\eta}(y,u(y),\D u(y))\cdot A(y) \Big\} \\
& \geq \ \H_{P}\big(x,u(x),\D u(x)\big) : \D A(x) \ +\ \H_{\eta}\big(x,u(x),\D u(x)\big) \cdot A(x) .
\end{split}
\]
If $\H_{P}\big(x,u(x),\D u(x)\big)\neq 0$, then by the definition of $\mathscr{L}\big( x, A(x), \X_x\big)$ we have
\[
\begin{split}
\underline{\D}r(0^+) &\geq \ \H_{P}\big(x,u(x),\D u(x)\big) : \D A(x) \ +\ \H_{\eta}\big(x,u(x),\D u(x)\big)\cdot A(x)
\\ 
&= \, -\, A(x)\cdot \Big(\mF_{\infty}^{\bot}\big( x,u(x),\D u(x),\X_x \big) \, - \,\H_{\eta}\big(x,u(x),\D u(x)\big)\Big) \\
&=\, -\, A(x) ^\top  [\![\H_{P}(x, u(x), \D u(x))]\!]^\bot \Big(  \mF_{\infty}^{\bot}\big(x,u(x), \D u(x), \X_x\big) \\
& \ \ \ \, - \, \H_{\eta}\big(x,\!u(x), \D u(x)\big) \Big)
\end{split}
\]
and hence $\underline{\D}r(0^+) \geq 0$ because $u$ is a $\mD$-solution on $\Om$. If $\H_{P}\big(x,u(x),\D u(x)\big)=0$, then again $\underline{\D}r(0^+) \geq 0$ because $A(x)=0$. In either cases, by inequality \eqref{2.1} we obtain $r(t) \geq 0$ for all $t\geq 0$ and hence
\[
\E_{\infty}(u,\,\mO)\, \leq\, \E_{\infty}(u+A,\,\mO), \ \ \ \forall\ \mO\Subset \Om, \  \forall\  A\in \mA^{\bot,\infty}_{\mO}(u).
\]
The theorem has been established.   \qed

\ms

\BPCOR \ref{corollary1}. If $u\in C^2(\Om,\R^N)$, then by Lemma \ref{lemma0} any diffuse hessian of $u$ satisfies $\mD^2u(x) = \de_{\D^2u(x)}$ for a.e.\ $x\in \Om$. By Remark \ref{Remark}, we may assume this happens for all $x\in \Om$. Therefore, the reduced support of $\mD^2u(x)$ is the singleton set $\{\de_{\D^2u(x)}\}$. Hence, for $ \smash{\mA^{\parallel,\infty}_{\mO}(u)}$, we have that any possible affine map $A$ satisfies $\D A\equiv \D\big( \xi \H\big(x,u(x),\D u(x)\big)\big)$ and $A(x)=0$. In the case of $\smash{\mA^{\bot,\infty}_{\mO}(u)}$, we have that any possible affine map  $A$ satisfies 
\[
A(x)^\top \H_{P}\big(x,u(x),\D u(x)\big)\, =\, 0\ , \ \ \ \D A\in \mathscr{L}\big( x , A(x) , \D^2u(x) \big), 
\]
which gives
\[
\begin{split}
\D A(x) : \H_{P}\big(x,u(x),\D u(x) & \big)\, =\, -A(x) \cdot \Big( \H_{PP}\big(x,u(x),\D u(x)\big):\D^2u(x) \, +
\\
& \hspace{-55pt} + \, \H_{P\eta}\big(x,u(x),\D u(x)\big): \D u(x) \, + \,  \H_{Px}\big(x,u(x),\D u(x)\big) : \mathrm{I} \Big)\\
& \ \ \,  =\,-A(x) \cdot \Div\big(\H_{P}\big(\cdot,u,\D u)\big)(x).
\end{split}
\]
As a consequence, the divergence $\mathrm{Div} \big( A^\top \H_{P}\big(\cdot,u,\D u\big) \big)(x)$ vanishes because
\[
\D A(x) : \H_{P}\big(x,u(x),\D u(x)\big) \, +\, A(x) \cdot \Div \big(\H_{P}(\cdot,u,\D u)\big)(x)\, =\, 0.
\]
The corollary has been established.        \qed

\ms
\ms

\noi \textbf{Acknowledgement.} N.K. would like to thank Jan Kristensen, Giles Shaw, Roger Moser, Tristan Pryer, Hussien Abugirda and Igor Velcic for inspiring scientific discussions on the topic of $L^\infty$ variational problems and of $\mD$-solutions to nonlinear PDE systems. He is also indebted to Gunnar Aronsson, Craig Evans, Juan Manfredi and Robert Jensen for their encouragement towards him. Both authors thank the anonymous referee for the careful reading of their manuscript and their constructive comments which led to substantial improvements of the content and the presentation of this paper.

\bibliographystyle{amsplain}

\end{document}